\documentclass[12pt]{amsart}
\usepackage{amssymb}
\usepackage{euscript}
\let\script\EuScript
\let\cal\mathcal

\let\epsilon\varepsilon

\makeatletter \@mparswitchfalse \makeatother

\newtheorem{theorem}{Theorem}
\newtheorem{lemma}[theorem]{Lemma}
\newtheorem{corollary}[theorem]{Corollary}
\newtheorem{proposition}[theorem]{Proposition}
\newtheorem{remark}[theorem]{Remark}
\newtheorem{remarks}[theorem]{Remarks}

\newtheorem{definition}[theorem]{Definition}

\def\eqref#1{(\ref{eq#1})}

\def\N{\mathbb N}

\def\R{\mathbb R}
\def\M{\cal{M}}
\let\<\langle

\def\T{\tau}

\numberwithin{equation}{section}
\numberwithin{theorem}{section}
\begin{document}

\let\\\cr

\let\union\bigcup
\let\inter\bigcap
\def\supp{\operatorname{supp}}
\def\sup{\operatorname{sup}}
\def\inf{\operatorname{inf}}

\def\log{\operatorname{log}}
\def\max{\operatorname{max}}
\def\min{\operatorname{min}}
\let\emptyset\varnothing
\def\exp{\operatorname{exp}}
\def\Sup{\operatornamewithlimits{Sup}}

\title[Non-commutative spaces]{ Sequences In
non-commutative $L^p$-spaces}
\author{Narcisse Randrianantoanina}
\thanks{Supported in part by NSF grant DMS-9703789}
\address{Department of Mathematics and Statistics, Miami University,
Oxford, Ohio 45056}
\email{randrin@muohio.edu}

\begin{abstract}
Let $\cal M$ be a semi-finite von Neumann algebra
equipped with a
distinguished faithful, normal, semi-finite trace
$\tau$.
We introduce the notion of equi-integrability in non-commutative spaces
and
show that if a rearrangement invariant quasi-Banach function space $E$
on the positive semi-axis is $\alpha$-convex with constant $1$ and
satisfies a non-trivial lower $q$-estimate with constant $1$, then the
corresponding
non-commutative space of measurable operators $E({\cal M}, \tau)$
 has the
following property:  every bounded sequence in
$E({\cal M}, \tau)$ has a subsequence
that splits into a $E$-equi-integrable sequence and a sequence with
pairwise disjoint projection supports.  This result extends the well known
Kadec-Pe\l czy\'nski subsequence decomposition
 for Banach lattices to  non-commutative
spaces. As applications, we prove that for $1\leq p <\infty$, every
subspace of $L^p(\cal M, \tau)$ either contains 
almost isometric copies of
$\ell^p$ or is strongly embedded in $L^p(\cal M, \tau)$.
\end{abstract}

\maketitle

\section{Introduction}
In \cite{KP}, Kadec and Pe\l czy\'nski proved that
if $1\leq p<\infty$ then every bounded sequence
$\{f_n\}_{n=1}^\infty$ in $L^p[0,1]$ has a
subsequence that can be decomposed  into
two extreme sequences $\{g_k\}_{k=1}^\infty$ and $\{h_k\}_{k=1}^\infty$, 
where the $h_k$'s are pairwise
disjoint and the $g_k$'s are $L_p$-equi-integrable that is
$\lim\limits_{m(A) \to 0}\sup\limits_{k}\Vert \chi_{A} g_{k}
 \Vert_p \rightarrow 0$ and $h_k \perp g_k$ for every $k \geq 1$.
This result was used to study different
structures of subspaces of $L^p[0,1]$.  Later, the same decomposition
property
was proved for larger classes of Banach function spaces (see \cite{FGJ}
 for Orlicz spaces with $\Delta_{2}$-condition and
$q$-concave lattices, \cite{JMST} for some symmetric spaces). There are
however Banach lattices with sequences for which the above decomposition
is not possible. Examples of reflexive, $p$-convex Banach lattices without
the subsequence decomposition can be found in a paper of Figiel
{\it et al} \cite{FGJ}.
Subsequently, Weis \cite{WE} characterized, in terms of uniform order
continuity conditions and ultrapowers, all Banach lattices where such
property is possible. For the case of rearrangement invariant function
spaces, the order-continuous spaces, in which  the above decomposition works, were
 fully characterized as those 
 that have the Fatou property (equivalently, those
that contains  no subspace isomorphic to
$c_{0}$).  In \cite{D3LRS}, a version of Kadec-Pe\l czynski decomposition
was considered for preduals of semi-finite von Neumann algebras.

It is the intention of the present paper to give an extension of the
Kadec-Pe\l czy\'nski decomposition stated above to the case of general
non-commutative
symmetric spaces of measurable operators.
 Let $\cal M$ be a von Neumann
algebra, equipped with a distinguished faithful, normal, semi-finite trace
$\tau$ and $E$ be a rearrangement invariant Banach function space
on $[0,1]$ or the half line $(0,\infty)$.  We define equi-integrability
in the non-commutative setting as generalizations of Akemann's
characterization of weak compactness on preduals of von Neumann algebras.
Using such notion, we provide an analogue of the Kadec-Pe\l czy\'nski
decomposition for non-commutative spaces. Namely, we proved
 that if $E$ is order continuous and satisifies the Fatou property then the
corresponding symmetric space of measurable operators
$E(\cal M,\tau)$ has the subsequence
splitting property.   Our approach
allows ones to consider more general spaces such as quasi-Banach
rearrangement invariant spaces that are $\alpha$-convex with constant $1$
and satisfy non trivial $q$-lower estimate with constant $1$.
In particular splitting of bounded sequences
is valid in non-commutative $L^p$-spaces for 
$0<p<\infty$. It should be noted that Sukochev
\cite{SU} obtain a similar result for the case of
finite von Neumann algebras. We also
remark that since $c_{0}$ fails the subsequence splitting property, our
result for the case where $E$ is a Banach space case 
is the best possible.

As application of the main result, we study the structure of
subspaces of $L^p(\M,\T)$ for $1<p<\infty$ which generalizes
the case $p=1$ treated in \cite{D3LRS}.

 We refer to
\cite{KR} and \cite{TAK} for general information concerning von Neumann
algebras as well as non-commutative integration, to \cite{SA} and
\cite{LT} for Banach lattice theory.

\section{Definitions and preliminary results}

Throughout, $H$ is a given Hilbert space and
${\cal M}\subset \cal{B}(H) $ denotes a semi-finite von Neumann algebra with a
distinguished normal,
faithful semi-finite trace $\tau$.  The identity in
$\cal M$ will be denoted
by $\bf{1}$ and $\M_p$ will stand for the set of all (self adjoint) projections
in $\M$.
 A closed and densely defined operator $a$ on $H$ is said to be
affiliated with $\cal M$ if $ua = au$ for all unitary  operator $u$ in the
commutant $\cal{M}'$ of $\cal M$.

A closed and densely defined operator $x$, affiliated with $\cal M$,
 is called
$\tau$-measurable if for every $\epsilon > 0$, there exists an orthogonal
projection $p \in \M$ such that $p(H)\subseteq
\text{dom}(x)$, $\tau ({\bf 1} -p)<\epsilon$ and $xp \in \cal{M}$.
The set of all $\tau$-measurable operators will be denoted by
$\widetilde{\cal M}$.  The set $\widetilde{\cal M}$ is a
${*}$-algebra with respect to the strong
sum, the strong product and the adjoint operation.
 Given a self-adjoint
operator $x$ in $\widetilde{\cal M}$, we denote by $e^{x}(\cdot)$ its spectral
measure.  Recall that $e^{\vert x \vert}(B) \in \cal{\cal M}$
 for all Borel sets $B
\subseteq \R$ and $x \in \widetilde{\cal M}$.
For fixed $x \in \widetilde{\cal M}$ and $t \geq 0$,
we define
$$\mu_{t}(x) = \inf \left\{s\geq 0 : \tau(e^{\vert x \vert}(s,
\infty))\leq t \right\}.$$
  The function  $\mu_{(.)}(x):[0, \infty) \rightarrow [0,
\infty]$ is called the generalized singular value function (or decreasing
rearrangement) of $x$. We note that $\mu_t(x) < \infty$ for every $t > 0$.
For a complete study of $\mu_{(.)}$, we refer the reader to \cite{FK}.
The topology defined by the metric on $\widetilde{\M}$ obtained by setting
$$d(x,y)=\inf\left\{t\geq 0:\ \mu_t(x-y)\leq t \right\}, 
\ \ \ \text{for}\ 
x, y \in \widetilde{\M},$$
is called the {\it measure topology}. It is well-known that a net
$(x_\alpha)_{\alpha \in I}$ in $\widetilde{\M}$  converge to 
$x \in \widetilde{\M}$ in measure  topology if and only if
for every $\epsilon >0$, $\delta >0$, there exists $\alpha_0 \in I$
such that whenever $\alpha \geq \alpha_0$, there exists a projection
$p \in \M_p$ such that
   $$\left\Vert (x_\alpha - x)p \right\Vert_{\M} <\epsilon \ \
\text{and}\ \ \T({\bf 1}-p)<\delta.$$
It was shown in \cite{N} that $(\widetilde{\M}, d)$ is 
a complete metric space.

Recall that if we consider $\cal{M} = L^{\infty}(\R^{+},m)$, where $m$ is the
Lebesgue measure on $\R^{+}$ then $\cal{M}$ is an abelian von Neumann algebra
acting on $L^2(\R^{+},m)$ via multiplication, with the trace being the
integration with respect to $m$, then
$\widetilde{\cal M} = L^{0}(\R^{+},m)$ (the usual
space of all measurable functions on $\R^{+}$) and the generalized
singular value $\mu(f)$ is precisely the decreasing rearrangement of the function
$\vert f\vert$ (usually denoted by $f^*$
in Banach lattice theory).

\begin{definition}
A symmetric quasi-Banach function space on $\R^+$ is a
quasi-Banach lattice
$E$ of measurable functions with the following properties:
\begin{itemize}
\item [(i)] $E$ is an order ideal in $L^{0}(\R^{+},m)$;
\item [(ii)] $E$ is rearrangemant invariant in the sense of
\cite{LT} (p~.114);
\item [(iii)] $E$ contains all finitely supported simple functions.
\end {itemize}
\end{definition}

Unless stated otherwise, $E$ will always denote a
 symmetric  quasi-Banach
function space on $\R^+$. We define the symmetric space of measurable
operators $E(\cal{ M}, \tau)$ by setting
$$E(\cal{M},\tau) := \left\{x \in \widetilde{\cal M} : \mu(x)\in E 
\right\}$$
and
$$\left\Vert x \right\Vert_{E(\cal{M}, \tau)} =
 \left\Vert \mu(x) \right\Vert_{E} \  \ \text{for \ all} \ x \in
E(\cal{M},\tau).$$

It is shown in \cite{X} (Lemma 4.1) that if $E$ is  $\alpha$-convex (for
some $0 < \alpha \leq 1$) with constant $1$, then $\Vert \cdot
\Vert_{E(\cal{M}, \tau)}$ is an $\alpha$-norm  that is
 for every $x, y \in E(\cal{M},\tau)$,
$$\left\Vert x +
y \right\Vert_{E(\cal{M}, \tau)}^\alpha \leq
 \left\Vert x \right\Vert_{E(\cal{M}, \tau)}^\alpha + \left\Vert y
\right\Vert_{E(\cal{M}, \tau)}^\alpha.$$
 Equipped with $\Vert \cdot \Vert_{E(\cal{M}, \tau)}$,
the space $E(\cal{M}, \tau)$ is
a $\alpha$-Banach space.  The space $E(\cal{M}, \tau)$ is
often referred to as the non-commutative analogue
of the function space $E$.  We remark that if
$0<p<\infty$ and $E = L^p
(\R^{+}, m)$  then
$E(\cal{M}, \tau)$ coincides with the usual
non-commutative $L^p$-space associated to the semi-finite von Neumann
algebra $\cal M$.  Also if
$E=L^\infty(\R^{+}, m)$, then 
$L^\infty(\M,\tau)$ is the von Neumann algebra
$\M$. We refer to \cite{DDP1}, \cite{DDP3} and \cite{X} for
some background on the space $E(\cal{M},\tau)$.

\begin{definition}
A quasi-Banach function space $E$ is said to satisfy a lower $q$-estimate
if there exists a positive constant $C>0$ such that for all finite
sequences $\{x_{n}\}$ of mututally disjoint elements in $E$,
$$\left(\sum \left\Vert x_n \right\Vert^q \right)^{\frac1{q}} \leq
C \left\Vert \sum x_n \right\Vert.$$
\end{definition}

The least such constant $C$ is called the constant of the lower
$q$-estimate.  Recall that if $E$ is a quasi-Banach function space and
$1<p<\infty$,
$$E^{(p)} = \left\{x \in L^0(\R^{+}, m); \vert x \vert^{p} \in E \right\}
 \  with \ \left\Vert x
\right\Vert_{E^{(p)}} = \left\Vert \  \vert x \vert^{p} \right\Vert_{E}^{\frac1{p}}.$$

We will  need the following known result. A proof can be found in
\cite{D3LS}.  
\begin{proposition}\label{ordercontinuity}
Assume that $E$ is order-continuous and $\alpha$-convex with constant $1$
 for some $0 < \alpha \leq 1$.
\begin{itemize}
\item [(i)] If $x \in E(\cal{M}, \tau)$ and $e \leq f$ are projections
in $\cal M$ then
$\Vert exe \Vert_{E(\cal{M}, \tau)} \leq \Vert fxf \Vert_{E(\cal{M}, \tau)}$;
\item [(ii)] If $x \in E(\cal{M}, \tau)$ and
$e_{\beta}\downarrow_\beta 0$ is a net of
 projections in $\cal M$ then
 $\Vert x e_\beta\Vert_{E(\cal{M}, \tau)} \downarrow_\beta 0$.
\end{itemize}
\end{proposition}

The following definition isolates the main topic of this paper.

\begin{definition}\label{equi-integrability}
Let  $E$ be a quasi-Banach function space on $\R^{+}$
and $K$ be a bounded subset of $E(\cal{M}, \tau)$.
We will say that $K$ is
{\bf $E$-equi-integrable} if $\lim\limits_{n \to \infty} \sup\limits_{x \in K}
\left\Vert e_n x e_n \right\Vert_{E(\cal{M},\tau)} = 0$
for every decreasing sequence $\{e_n\}_{n=1}^\infty$ of projections with
$e_n \downarrow_n  0$.
\end{definition}

We remark, from Proposition~\ref{ordercontinuity}, that since 
$\{e_n\}_{n=1}^\infty$ is decreasing, so is the sequence
$\left\{\sup\limits_{x \in K}
\left\Vert e_n x e_n \right\Vert_{E(\cal{M},\tau)}\right\}_{n=1}^\infty$
and therefore the limit in the definition above always exists.
This definition was motivated by the commutative case on one hand
and the characterization of weakly compact subsets of
$L^{1}(\cal{M}, \tau)$ by Akemann \cite{AK} (see also
\cite{TAK} p.150) on the other. Using
this terminology, Akemann's characterization can be stated as in the
commutative case:  relatively weakly compact subsets
of $L^1(\cal{M}, \tau)$ are exactly the
equi-integrable sets. Such characterization is not valid 
in general. For $1<p<\infty$, any set of normalized disjoint
sequence cannot be $L^p$-integrable but since $L^p$ is 
reflexive, such set is relatively weakly compact.

On the next proposition, we will show that for the general
case, one implication always holds.

\begin{proposition}\label{compact}
Assume that $E$ is an order-continuous symmetric 
 Banach function space and $K$ is a $E$-equi-integrable 
set in $E(\cal{M},\tau)$ then $K$ is relatively weakly compact.
\end{proposition}

The proposition will be proved in several steps. Recall that if $E$ is a
symmetric Banach function space, then $E(\cal{M},\tau)$ is a subset of
$L^1(\cal{M} , \tau) + \cal{M}$ and therefore if $p$ is a projection in
$L^1(\cal{M}, \tau ) \cap \cal{M}$  and $K$ is a subset of
$E(\cal{M},\tau)$, then $pK$ and $Kp$ are subsets of $L^1(\cal{M},\tau)$.
\begin{lemma}\label{finitecompact1}
Let $p$ be a projection in $L^1(\cal{M},\tau) \cap \cal{M}$ and $K$
be a $E$-equi-integrable subset of $E(\cal{M},\tau)$. The sets
$pKp$ and $pK(1-p)$ are relatively weakly compact in $L^1(\cal{M},\tau)$.
\end{lemma}
To see this lemma, it is enough to check that these sets are
$L^1$-equi-integrable.
Let $T: E(\cal{M},\tau) \longrightarrow L^1(\cal{M},\tau)$ be the linear
map defined by $x \to Tx= pxp$. This map is well-defined and one can
deduce from the closed graph theorem that it is bounded. Let 
$\{e_n\}_{n=1}^\infty$ be a
sequence of  projections with $e_n \downarrow_n 0$.
For each $n\geq 1$, set $f_n$ to be the right support projection
of $e_np$. By the definition of support projections, $f_n \leq p$. So
$\{f_n\}_{n=1}^\infty$ is a sequence of finite projections. We also note that
(see for instance the proof of \cite[Proposition~1.6 p.292]{TAK}),
$$f_n =  e_n \vee ({\bf 1}-p) - ({\bf 1}-p)$$
and by Kaplansky formula (see for instance \cite[Theorem~6.1.6 p.403]{KR}),
 $$f_n \sim  e_n - e_n \wedge ({\bf 1}-p).$$
 Since $\tau(f_n) =\tau(e_n - e_n \wedge({\bf 1}-p)) \leq \tau(p)$ and
 $\{e_n - e_n \wedge({\bf 1}-p)\}_{n=1}^\infty$ converges to zero, 
$\{f_n\}_{n=1}^\infty$ converges to
 zero. Now since the $f_n$'s are finite projections, we conclude that
 if $g_n= \wedge_{k \geq n}f_k $, then $\{g_n\}_{n=1}^\infty$ converges to zero. Therefore, for every
 $x \in K$,
 \begin{align*}
  \left\Vert e_n pxp e_n \right\Vert_1 
&=\left\Vert e_n p (g_n x g_n) pe_n \right\Vert_1 \cr
  &\leq \left\Vert  p(g_n x g_n)p \right\Vert_1 \cr
  &\leq ||T||  \cdot \left\Vert g_n x g_n \right\Vert_{E(\cal{M}, \tau)}.
\end{align*}
Since $K$ is $E$-equi-integrable, one obtain that
$$\lim_{n\to \infty}\sup_{y \in pKp} \left\Vert e_n y e_n \right\Vert_1
\leq ||T|| \cdot \lim_{n \to \infty}\sup_{x\in K}\left\Vert g_n x g_n
\right\Vert_{E(\cal{M},\tau)}=0 $$
which concludes that $pKp$ is relatively weakly compact in
$L^1(\cal{M},\tau)$.

For $pK({\bf 1}-p)$, set $S: E({\cal{M}},\tau) \to L^1(\M,\T)$ 
be the map
defined by $x \to Sx=px({\bf 1}-p)$. As above, $S$ is bounded. Let 
$\{e_n\}_{n=1}^\infty$ and
$\{g_n\}_{n=1}^\infty$ be sequences of projections as discribed above. For each
$ n\geq 1$, let $s_n$ be  the left support projection of 
$({\bf 1}-p)e_n$. Then 
$s_n =e_n \vee p - p $
for every $n \geq 1$ and the sequence $\{s_n\}_{n=1}^\infty$ is decreasing. It is
 claimed
that $s_n \downarrow_n 0$.

For this, it is enough to check that $e_n \vee p \downarrow_n p$.
In fact, $e_n \vee p - e_n \sim p - e_n \wedge p$ and the sequence
defined by the right hand side of the equivalence
converges to $p$ which implies  that
$$\lim_{n\to \infty} \tau(e_n \vee p - e_n)= \lim_{n \to \infty}
\tau(p- e_n\wedge p) =\tau(p)$$
which gives
$$\lim_{n \to \infty} \tau(e_n \vee p - p - e_n)=0.$$
But since 
$(e_n \vee p - p - e_n)^2 =(e_n \vee p - p - e_n) + e_np + pe_n$,
we can conclude that
$$\lim_{n \to \infty} \Vert e_n \vee p - p - e_n \Vert_2 =0.$$
From this, we get (by passing to a subsequence if necessary) that
$\{e_n \vee p - p - e_n\}_{n=1}^\infty$ converges to zero in measure. Similarly,
$\{e_n \vee p - p -pe_n\}_{n=1}^\infty$ converges to zero in measure so
$e_n \vee p \downarrow_n p$  hence $s_n \downarrow_n 0$.

To conclude the proof of Lemma~\ref{finitecompact1}, note that $g_n \perp s_n$ so
$g_n \vee s_n = g_n + s_n $. In particular,
$g_n \vee s_n \downarrow_n 0$ and we get that
\begin{equation*}
\begin{split}
\lim_{n\to \infty}\sup_{y \in pK({\bf 1}-p)} 
\left\Vert e_n y e_n
\right\Vert_1
&=\lim_{n\to \infty}\sup_{x \in K} \left\Vert e_n p x 
({\bf 1}-p)
e_n \right\Vert_1\\
&=\lim_{n\to \infty}\sup_{x \in K} \left\Vert e_n p
(g_n \vee s_n) x (g_n \vee s_n)({\bf 1}-p)  e_n \right\Vert_1 \\
&\leq\lim_{n\to \infty}\sup_{x \in K} \left\Vert  p
(g_n \vee s_n) x (g_n \vee s_n)({\bf 1}-p)   \right\Vert_1 \\
&\leq \|S\|\cdot \lim_{n\to \infty}\sup_{x \in K} \left\Vert
 (g_n \vee s_n) x (g_n \vee s_n)   \right\Vert_{E(\M,\tau)} =0.
\end{split}
\end{equation*}

The proof  is complete.

\begin{lemma}\label{finitecompact2}
Let $p$ and $K$ be as in Lemma~\ref{finitecompact1}. Then $pK$ is relatively weakly
compact in $E(\M,\T)$.
\end{lemma}

Note first that $pK$ is $E$-equi-integrable. This can be seen  by applying
the series of argument used in Lemma~\ref{finitecompact1}, considering
 the operators $T$ and $S$
as maps from  $E(\M,\T)$ into $E(\M,\T)$.

Let $\{px_n\}_{n=1}^\infty$ be a bounded sequence in $pK$. 
From Lemma~ref{finitecompact1},
 we can assume
that $\{px_n\}_{n=1}^\infty$ is weakly convergent in $L^1(\M,\T)$. Fix $\varphi \in
E^*(\M, \T)_{+}$  and let $\varphi=\int_{0}^\infty t \ de_t$ be its
spectral decomposition. For each $k \geq 1$, set
$q_k:=e^{\varphi}( (0,k))$.
 We remark that $\varphi q_k = q_k \varphi \in \M$.
For $m, \ n \in \N$,
\begin{align*}
\langle \varphi, px_n -px_m \rangle
&=\langle \varphi -\varphi q_k, px_n -px_m \rangle +
 \langle \varphi q_k, px_n -px_m \rangle \\
&=\langle \varphi ({\bf 1}-q_k), px_n -px_m \rangle +
 \langle \varphi q_k, px_n -px_m \rangle \\
&=\langle ({\bf 1}-q_k)\varphi ({\bf 1}-q_k), px_n -px_m \rangle +
 \langle \varphi q_k, px_n -px_m \rangle.
\end{align*}
This gives
 $$ \left|\langle \varphi, px_n -px_m \rangle\right| \leq
 \left|\T\left( \varphi (1-q_k)(px_n -px_m)(1-q_k) \right) \right| +
 \left|\langle \varphi q_k, px_n -px_m \rangle \right|.$$
 Since $\varphi q_k$ belongs to $\M$,
 $$\limsup\limits_{n,m \to \infty}\left|\langle \varphi q_k, px_n -px_m
\rangle\right|
 \leq 2 \left\Vert\varphi \right\Vert_{E^*(\M,\T)}
\cdot \sup_{a \in K}\left\Vert ({\bf 1}-q_k)pa({\bf 1}-q_k)
 \right\Vert_{E(\M,\T)}.$$
 Now since ${\bf 1}-q_k \downarrow_k 0$  and $pK$ is $E$-equi-integrable, we
 obtain that
  $$\lim_{n,m \to \infty}\langle \varphi q_k, px_n -px_m
  \rangle=0.$$
  The lemma is proved.

\smallskip

To complete the proof of Proposition~\ref{compact}, let 
$\{p_k\}_{k=1}^\infty$ be a sequence
of projections that increases to ${\bf 1}$ and $\T(p_k)< \infty$ and fix
$\epsilon > 0$. Choose $k_0 \geq 1$ such that
$$\sup_{a \in K} \left\Vert ({\bf 1}-p_{k_0})a({\bf 1}-p_{k_0})\right\Vert_{E(\M,\T)} \leq
\epsilon.$$
We have $K=p_{k_0}K + ({\bf 1}-p_{k_0})Kp_{k_0} + 
({\bf 1}-p_{k_0})K({\bf 1}-p_{k_0})$ which
implies that
 $$K \subset p_{k_0}K + ({\bf 1}-p_{k_0})Kp_{k_0} +  \epsilon B_{E(\M,\T)}$$
  where
 $B_{E(\M,\T)}$ denotes the closed unit ball of $E(\M,\T)$.
 From Lemma~\ref{finitecompact2}, the sets 
$p_{k_0}K$ and $({\bf 1}-p_{k_0})Kp_{k_0}$ are relatively
 weakly compact which concludes that $K$ is relatively weakly compact.
 The proof is complete.
\qed

\begin{remark}
If $\tau({\bf 1}) <\infty$, the proof above can be considerably shortened.
In this case, $E(\M,\T) \subset L^1(\M,\T)$ so if $K$ is $E$-equi-integrable,
then it is relatively weakly compact in $L^1(\M,\T)$ and on can argue as in 
the last part of Lemma~\ref{finitecompact2}
 to conclude that $K$ is relatively weakly
compact in $E(\M,\T)$.
\end{remark}

The following proposition should be compared with
\cite{CS2}(Theorem~5.1 and Theorem~5.2). It generalizes well known property
of equi-integrable sets in function spaces to the non-commutative settings. 
\begin{proposition}\label{measure}
Let $E$ be a symmetric quasi-Banach space function space and $K$ be a 
$E$-equi-integrable subset of $E(\M,\T)$.
For each sequence $\{x_n\}_{n=1}^\infty$ in $K$ and $x \in \overline{K}$,
the following are equivalent:
\begin{itemize}
\item[(a)] $\lim\limits_{n \to \infty} \left\Vert x_n -x
\right\Vert_{E(\M,\T)} =0$;
\item[(b)] $\{x_n\}_{n=1}^\infty$ converges to $x$ in measure (as 
$n \to \infty$).
\end{itemize}
\end{proposition}

\begin{proof}
The implication $(a) \Rightarrow (b)$ is trivial.
For $(b) \Rightarrow (a)$, we will assume that $x=0$.
Recall that there exists $0 < \alpha \leq 1$, such that
$\M \cap L^\alpha(\M,\T) \subset E(\M,\T) \subset 
\M + L^\alpha(\M,\T)$ with
$\Vert x \Vert_{\M +L^\alpha(\M,\T)}
\leq \Vert x \Vert_{E(\M,\T)} \leq
2\Vert x \Vert_{\M \cap L^\alpha(\M,\T)}$ for every
$x \in \M \cap L^\alpha(\M,\T)$.
The proposition will be proved by
showing the following lemma:
\begin{lemma}
For every $p \in \M_p $ with $\T(p)<\infty$,
$\lim\limits_{n \to \infty}
\left\Vert x_np \right\Vert_{E(\M,\T)} =
\lim\limits_{n \to \infty}
\left\Vert px_n \right\Vert_{E(\M,\T)} =0$.
\end{lemma}
To see this lemma, fix $\epsilon >0$ and let
$C=\max\{1,\tau(p)\}$. Since $K$ is
equi-integrable, there exists $\delta>0$ such that
whenever $q \in \M_p$ satisfies $\tau(q)<\delta$, then
for every $n \in \N$,
$\left\Vert qx_nq \right\Vert_{E(\M,\T)} \leq \epsilon/(2)^{1/\alpha}$.
Since both $\{x_n\}_{n=1}^\infty$ and $\{x_n^*\}_{n=1}^\infty$
converge to zero
in measure, one can choose $n_0 \geq 1$ such that for each
$n \geq n_0$, there exists a projection $p_n \in \M_p$
with $\T({\bf 1}-p_n) <\delta$,
$$\left\Vert x_n p_n \right\Vert_{\M} <
\frac{\epsilon}{2[4C]^\frac{1}{\alpha}}$$
and
$$\left\Vert x_n^*p_n \right\Vert_{\M} <
\frac{\epsilon}{2[4C]^\frac{1}{\alpha}}.$$
For $n \geq n_0$,
\begin{align*}
 \left\Vert x_np \right\Vert_{E(\M,\T)}^\alpha
&\leq \left\Vert x_n p_n p \right\Vert_{E(\M,\T)}^\alpha
+ \left\Vert p_n x_n ({\bf 1}-p_n) p
\right\Vert_{E(\M,\T)}^\alpha
+ \left\Vert ({\bf 1}-p_n) x_n ({\bf 1}-p_n) p
\right\Vert_{E(\M,\T)}^\alpha \\
&\leq 2^\alpha \max\left\{\left\Vert x_n p_n
\right\Vert_{\M }^\alpha,
    \left\Vert x_n p_n p \right\Vert_{L^\alpha(\M,\T)}^\alpha
\right\} \\
&\ +
2^\alpha
 \max\left\{\left\Vert x_n^*p_n \right\Vert_{\M }^\alpha,
    \left\Vert p({\bf 1}-p_n)x_n^*p_n \right\Vert_{L^\alpha(\M,\T)}^\alpha
\right\}
 +
\left\Vert ({\bf 1}-p_n) x_n ({\bf 1}-p_n)
\right\Vert_{E(\M,\T)}^\alpha \\
&\leq 2. 2^\alpha \max\left\{ \epsilon^\alpha/ 2^\alpha 4C,
 (\epsilon^\alpha/2^\alpha 4C) \tau(p)  \right\} +
\epsilon^{\alpha}/2 \\
&\leq \epsilon^\alpha.
\end{align*}
A similar estimate works for $\{x_n^*p\}_{n=1}^\infty$.
The lemma is verified.

To complete the proof of Proposition~\ref{measure}, choose a mutually
disjoint family
$\{e_i\}_{i\in I}$ of
 projections in $\M$ with
$\sum_{i \in I} e_i={\bf 1}$ for the strong operator topology
 and $\tau(e_i)<\infty$ for all $i \in I$. Using a similar
argument as in \cite{X}, one can get an at most countable
subset $\{e_k\}_{k=1}^\infty$ of $\{e_i\}_{i \in I}$ such that for each
$e_i$ outside of $\{e_k\}_{k=1}^\infty$,
$e_ix_n =x_ne_i=0$ for every $n \in \N$.
Let $e=\sum_{k \in \N} e_k$. Replacing $\M$ by $e{\M}e$ and
$\tau$ by its restriction on $e{\M}e$, we may assume that
$e={\bf 1}$. Let $p_n=\sum_{k\geq n}e_k$. It is clear that
$p_n \downarrow_n 0$ and $\tau({\bf 1}-p_n) <\infty$ for
every $n \in \N$.  Fix $\epsilon>0$ and choose $n_0 \geq 1$
such that
 $$\sup_{n \in \N}\left\Vert p_{n_0} x_n p_{n_0}
\right\Vert_{E(\M,\T)} \leq \epsilon.$$
We get that
$$\limsup_{n\to \infty}
\left\Vert x_n \right\Vert_{E(\M,\T)}^\alpha
\leq \lim_{n\to \infty}
\left\Vert x_n({\bf 1}-p_{n_0}) \right\Vert_{E(\M,\T)}^\alpha +
\lim_{n\to \infty}
\left\Vert x_n^*({\bf 1}-p_{n_0} \right\Vert_{E(\M,\T)}^\alpha +
\epsilon =\epsilon$$
and since $\epsilon$ is arbitrary,
the proof   is complete.
\end{proof}

\bigskip

For the rest of this section, we collect some results that will be useful in
the later sections of the paper.

\smallskip

The inequality given below can be viewed as the analogue of
the well-known fact on normal functional on von Neumann algebra,
$|\varphi(a)|^2 \leq \Vert \varphi \Vert \cdot \vert \varphi
\vert(aa^*)$ whenever $a \in \M$ and $\varphi \in \M_*$
\cite[Proposition 4.6 p.~146]{TAK}, to the  general 
case of symmetric spaces
of measurable operators.

\begin{proposition}\label{inequality}
Let $x \in E(\cal{M},\tau)$ and $y \in \cal{M}$ then
$$\left\Vert xy \right\Vert_{E(\cal{M},\tau)} \leq \left\Vert \
\vert x \vert y \right\Vert_{E(\cal{M}, \tau)}
\leq \left\Vert x \right\Vert_{E(\cal{M},\tau)}^{\frac1{2}}
\cdot \left\Vert y^{*}\vert x \vert y
\right\Vert_{E(\cal{M},\tau)}^{\frac1{2}}.$$
\end{proposition}

\begin{proof}
Let $x = u \vert x \vert$ be the polar decomposition of $x$.  Then
$\left\Vert xy
\right\Vert_{E(\cal{M},\T)} = \left\Vert u \vert x \vert y 
\right\Vert_{E(\cal{M},\T)}
\leq \left\Vert u
\right\Vert_{\infty} \cdot \left\Vert \  \vert x \vert y 
\right\Vert_{E(\cal{M},\T)}$.
 Also $\left\Vert \  \vert
x \vert y \right\Vert_{E(\cal{M},\T)} = \left\Vert \  \vert x \vert^{\frac1{2}} \vert x
\vert^{\frac1{2}} y \right\Vert_{E(\cal{M},\T)}$
and using H\"older's inequality,
\begin{align*}
\left\Vert \  \vert x \vert y \right\Vert_{E(\cal{M},\T)}
&\leq \left\Vert \  \vert x \vert^{\frac1{2}}
\right\Vert_{E^{(2)}{(\cal{M},\T)}} \cdot \left\Vert \ 
\vert x \vert^{\frac1{2}} y
\right\Vert_{E^{(2)}{(\cal{M},\T)}}\cr
&= \left\Vert x \right\Vert_{E(\cal{M},\T)}^{\frac{1}{2}}
  \cdot \left\Vert y^* \vert x \vert y
\right\Vert_{E(\cal{M},\T)}^{\frac{1}{2}}.
\end{align*}
\end{proof}

\begin{remark}
Let  $K$ be a bounded subset of $E(\cal{M},\tau)$. If 
we set
$|K|:=\left\{ |a|:\ a \in K \right\}$, then
 it is clear from Proposition~\ref{inequality} that if 
 $\vert K \vert$ is
$E$-equi-integrable then for every decreasing projections
$e_n{\downarrow}_n 0$,
$\lim\limits_{n \to \infty} \sup\limits_{x \in K} 
\left\Vert x e_n
\right\Vert_{E(\M,\tau)}
= \lim\limits_{n \to \infty} \sup\limits_{x \in K} 
\left\Vert e_n x \right\Vert_{E(\M,\tau)}= 0$.
In particular if $|K|$ is $E$-equi-integrable then so is $K$.
\end{remark}

\begin{proposition}\label{equi-integrability2}  Assume that $E$ is $\alpha$-convex with
constant 1 for some $0 < \alpha \leq 1$.
Let $\{p_n\}_{n=1}^\infty$ be a sequence of decreasing projections in $\cal M$
and $K$ be a bounded subset of $E(\cal{M},\tau)$ such that:
\begin{itemize}
\item [(i)] $p_n \downarrow_n 0$;
\item [(ii)] For each $n \geq 1$, the sets
$({\bf 1} - p_n)K$ and $\vert K({\bf 1}-p_n) \vert$ are
$E$-equi-integrable.
\end{itemize}
Then $K$ is $E$-equi-integrable if and only if
$\lim\limits_{n \to \infty}
\sup\limits_{a \in K} \left\Vert p_n a p_n \right\Vert_{E(\cal{M},\tau)} = 0.$
\end{proposition}

\begin{proof}
We will show  the non trivial implication.
 Fix $f_k \downarrow_k 0$, a sequence
 in $\cal{M}_p$.   We need to show that
$\lim\limits_{k \to \infty} \sup\limits_{a \in K} \left\Vert f_k a f_k
\right\Vert_{E(\cal{M},\tau)} = 0$.

We will assume without loss of generality 
that $K$ is a subset of the unit ball of $E(\M,\T)$.
For every
$a \in K$,
\begin{align*}
f_k a f_k &= f_k({\bf 1}-p_n)af_k + f_k p_n a f_k \cr
&=f_k ({\bf 1}-p_n)a f_k + f_k p_n a ({\bf 1}-p_n)f_k + f_kp_n a p_n f_k.
\end{align*}
Since $E(\cal{M},\tau)$ is $\alpha$-convex, we get:
\begin{align*}
\left\Vert f_k a f_k \right\Vert^{\alpha}_{E(\cal{M},\tau)}
&\leq \left\Vert f_k ({\bf 1}-p_n)a f_k
\right\Vert^{\alpha}_{E(\cal{M},\tau)}
+ \left\Vert f_k p_n a({\bf 1}-p_n)f_k \right\Vert^{\alpha}_{E(\cal{M},\tau)}
 + \left\Vert f_k p_n a p_n f_k \right\Vert^{\alpha}_{E(\cal{M},\tau)} \cr
 &\leq \left\Vert f_k ({\bf 1}-p_n) a f_k \right\Vert^{\alpha}_{E(\cal{M},\tau)} +
 \left\Vert a ({\bf 1}-p_n) f_k
\right\Vert^{\alpha}_{E(\cal{M},\tau)}
 + \left\Vert p_n a p_n \right\Vert^{\alpha}_{E(\cal{M},\tau)}.
\end{align*}
Using Proposition~\ref{inequality} on the second term, we have
\begin{align*}
\left\Vert f_k a f_k \right\Vert^{\alpha}_{E(\cal{M},\tau)}
&\leq \left\Vert f_k ({\bf 1}-p_n)a f_k
\right\Vert^{\alpha}_{E(\cal{M},\tau)} +
\left\Vert a ({\bf 1}-p_n) \right\Vert^{\frac\alpha{2}}_{E(\cal{M},\tau)}
 \cdot \left\Vert
f_k \vert a({\bf 1}-p_n) \right\vert f_k \Vert^{\frac\alpha{2}}_{E(\cal{M},\tau)} \cr
&+ \left\Vert p_n a p_n \right\Vert^{\alpha}_{E(\cal{M},\tau)}.
\end{align*}
Let $\epsilon > 0$, choose $n_0$ large enough so that
$\sup\limits_{a \in
K} \left\Vert p_{n_{0}}a p_{n_{0}} \right\Vert_{E(\cal{M},\tau)}
<\epsilon$. We conclude that
\begin{align*}
\lim\limits_{k \to \infty} \sup\limits_{a \in K}
\left\Vert f_k a f_k
\right\Vert^{\alpha}_{E(\cal{M},\tau)}
&\leq \lim\limits_{k \to \infty} \sup\limits_{a \in
K} \left\Vert f_k ({\bf 1}-p_{n_0})a f_k \right\Vert^{\alpha}_{E(\cal{M},\tau)}\cr
&+ \lim\limits_{k \to \infty} \sup\limits_{a \in K} \left\Vert f_k \vert
a({\bf 1}-p_{n_0}) \vert f_k 
\right\Vert^{\frac\alpha{2}}_{E(\cal{M},\tau)}
 + \epsilon^\alpha.
\end{align*}
By (ii), the first two terms converge to zero so $\lim\limits_{k \to \infty}
\sup\limits_{a \in K} \Vert f_k a f_k \Vert_{E(\cal{M},\tau)} \leq \epsilon$
and since $\epsilon$ is arbitrary, the proof is complete.
\end{proof}

\smallskip

The next  proposition  
 can be found in \cite{D3LS}
(Proposition~2.5).

\begin{proposition}\label{lower-estimate}
Assume that $E$ is $\alpha$-convex with constant $1$ for  some
$0 < \alpha \leq
1$ and satisfies a lower $q$-estimate with constant $1$
for some finite $q \geq \alpha$.  If $k =2q/\alpha$, then for all
$y \in E(\cal{M}, \tau)$, for all projections $e, f \in \cal{M}$ with $e + f = 1$ and $\tau(e) <
\infty$, it follows that
$$
\Vert e y e \Vert^k_{E(\cal{M},\tau)} +
\Vert e y f \Vert^k_{E(\cal{M},\tau)} +
\Vert f y e\Vert^k_{E(\cal{M},\tau)} +
 \Vert f y f \Vert^k_{E(\cal{M},\tau)}
\leq \Vert y \Vert^k_{E(\cal{M}, \tau)}.
$$
\end{proposition}

\smallskip

The proof of the next lemma is just a notational adjustment of the proof of
a lemma  from  the $L^1$-case \cite{D3LRS} so we will leave
 the details to the interested readers.

\smallskip

\begin{lemma}\label{subsequence}
Let $\{x_k\}_{k=1}^\infty$ be a bounded sequence in $E(\cal{M},\tau)$ 
and $\{e_n\}_{n=1}^\infty$ be a
decreasing sequence of projections in $ \M$
such that $e_n\downarrow_n 0$.
Assume that $\lim\limits_{n \to \infty} \sup\limits_{k \in \N} \Vert e_n
x_k e_n \Vert_{E(\cal{M},\tau)} =\gamma > 0$ then there exists a subsequence
 $\{x_{k_n}\}_{n=1}^\infty$ so
that $\lim\limits_{n \to \infty} \Vert e_n x_{k_n} e_n \Vert_{E(\cal{M},\tau)} = \gamma$.
\end{lemma}

\section{Kadec-Pe\l czy\'nski theorem for symmetric spaces}

The main result of the present  article is the following theorem.

\begin{theorem}\label{main}
Let $E$ be an order continuous symmetric quasi-Banach function space in $\R^+$ that
is $\alpha$-convex with constant $1$ for some $0 < \alpha \leq 1$ and
suppose that $E$ satisfies a lower $q$-estimate with constant $1$ for some
$q \geq \alpha$.

Let $\{x_n\}_{n=1}^\infty$ be a bounded sequence in $E(\cal{M},\tau)$
 then there exists a
subsequence $\{x_{n_k}\}_{k=1}^\infty$ of $\{x_n\}_{n=1}^\infty$, 
bounded sequences $\{y_k\}_{k=1}^\infty$ and
$\{z_k\}_{k=1}^\infty$ in $E(\cal{M},\tau)$ and a decreasing sequence of
projections
 $p_k \downarrow_n 0$ in $\cal M$ such that:

\begin{itemize}
\item [(i)] $x_{n_k} = y_k + z_k$ for all $k \geq 1$;
\item [(ii)] $\{y_k:\ k \geq 1\}$ is $E$-equi-integrable and $p_ky_kp_k = 0$
for all $k \geq 1$;
\item [(iii)] $\{z_k\}_{k=1}^\infty$ is such that $p_kz_kp_k = z_k$ for all $k \geq
1$.
\end{itemize}
\end{theorem}

The proof will be divided into several steps.
Without loss of generality, we will assume that 
the sequence $\{x_k\}_{k=1}^\infty$ is a subset of the unit ball of $E(\M,\T)$. 
Since we are dealing with sequences, we can and do assume 
without loss of generality that
$\cal M$ is countably decomposable (see \cite{X} and the proof of
Proposition~\ref{measure} above 
for the details of such reduction).

Set $\script{D}_1 : = \left\{\{e_n\}_{n=1}^\infty \subset \cal{M}_p:\  e_n \downarrow_n 0
\ \text{and}\  \tau(e_1)< \infty \right\}$ 
and consider
 $$\delta: = \sup\left\{\lim\limits_{n \to \infty} \sup\limits_{k \in \N}
\left\Vert e_n \vert x_k \vert e_n \right\Vert_{E(\cal{M},\tau)}:\ 
\{e_n\}_{n=1}^\infty \in \script{D}_1 \right\}.$$

As in \cite{D3LRS}, one can show that $\delta$ is attained and using
Lemma~\ref{subsequence},
one can choose a subsequence of $\{x_{n}\}_{n=1}^\infty$ (which we will
denote again by $\{x_n\}_{n=1}^\infty$) and $\{e_n\}_{n=1}^\infty \in \script{D}_1$ 
such that
\begin{equation}\label{delta}
\delta =\lim\limits_{n \to \infty} \left\Vert e_n \vert x_n \vert e_n
\right\Vert_{E(\cal{M},\tau)}
\end{equation}
 and
\begin{equation}\label{delta*}
\delta^*=\lim\limits_{n \to \infty} \left\Vert e_n \vert x^{*}_{n} \vert e_n
\right\Vert_{E(\cal{M},\tau)}
= \sup \left\{\lim\limits_{n \to \infty} \sup\limits_{k} 
\left\Vert q_n \vert
x^{*}_{k} \vert q_n \right\Vert_{E(\cal{M},\tau)}; 
\{q_{n}\}_{n=1}^\infty \in \script{D}_1 \right\}.
\end{equation}

For each $n \geq 1$, set
$v_n: = x_n - e_n x_n e_n$ and let
$V: = \left\{v_n:\  n \geq 1 \right\}$.

\begin{lemma}
There exists a sequence of projections $\{g_{n}\}_{n=1}^\infty$ in $\cal{M}$ with:
\begin{itemize}
\item [(1)] For every $n \geq 1$, $g_n \leq {\bf 1} - e_n$;
\item [(2)]  $\tau(g_n)<\infty$, in particular
$g_n$ is a finite projection;
\item [(3)] $g_n \uparrow^n \bf{1}$.
\end{itemize}
\end{lemma}

\begin{proof}
The lemma can be obtained inductively.
Since $\cal M$ is countably decomposable, there exists
  $\varphi_0$  a faithful normal
state in $\cal{M}_*$.  Since ${\bf 1} - e_n$ is a
semifinite projection, there exists a sequence of
projections $\{g_j^{(n)}\}_{j=1}^\infty$ with $\tau (g_j^{(n)}) < \infty$ for every $j \geq
1$ and $g_j^{(n)} \uparrow^j {{\bf 1}-e_n}$.  One can choose $j_n \geq 1$ such that
$\varphi_0 (1-e_n) - \varphi_0(g^{(n)}_{j_n})< 1/n$.
Set
$$\begin{cases}
g_n: = g^{(1)}_{j_1} &\text{for } n=1 \cr
 g_n = g^{(n)}_{j_n} \vee g_{n-1} &\text{for } n>1.
\end{cases}$$
 It is easy to verify  that $\{g_n\}_{n=1}^\infty$ satisfies the requirements of 
the lemma.
\end{proof}

\smallskip

For each $n \geq 1$, let
$p_n = {\bf 1} - g_n$.  Clearly $p_n \downarrow_{n} 0$, 
${\bf 1}-p_{n}$ is a finite
projection and $p_n \geq e_n$ for each $n \geq 1$.

\begin{lemma}\label{finite-equi-integrable}
For each $n \geq 1$, the sets 
$\left\vert V ({\bf 1}-p_n) \right\vert$ and 
$\left\vert({\bf 1}-p_n)V\right\vert$ are
$E$-equi-integrable.
\end{lemma}

\begin{proof}
Let us show that for every $n \geq 1$, 
$\left\vert V (1-p_{n}) \right\vert$ is an $E$-equi-integrable set.
Assume that there exists $ k_0 \geq 1$ such that $\vert V
({\bf 1}-p_{k_0}) \vert$ is not $E$-equi-integrable.
 By definition, there exists a
decreasing sequence of projections $q_n \downarrow_n 0$
such that
$\lim\limits_{n \to \infty} \sup\limits_{a \in \vert V ({\bf 1}-p_{k_0})\vert}
\left\Vert q_n\ a\ q_n \right\Vert_{E(\cal{M},\tau)} > 0$, that is
$$\lim_{n \to \infty} \sup_{m \in \N}
\left\Vert q_n \vert v_m
({\bf 1}-p_{k_0}) \vert q_n \right\Vert_{E(\cal{M},\tau)} > 0.$$
Choose a strictly increasing sequence $\{m_n\}_{n=1}^\infty$ of
  $\N$ such that
$$\lim_{n \to \infty} \left\Vert q_n \vert v_{m_n}(1-p_{k_0})
\vert q_n \right\Vert_{E(\cal{M},\tau)} > 0.
$$
Let $u_{n,k_0}$ be a bounded operator such that
$\vert v_{m_n} ({\bf 1}-p_{k_0})
\vert = u_{n,k_0} v_{m_n} ({\bf 1}-p_{k_0})$.
We get that
\begin{align*}
\left\Vert q_n \vert v_{m_n} ({\bf 1}-p_{k_0}) \vert q_n
\right\Vert_{E(\cal{M},\tau)} &=
\left\Vert q_n u_{n,k_0} v_{m_n} ({\bf 1}-p_{k_0}) q_n
\right\Vert_{E(\cal{M},\tau)} \cr
&= \left\Vert q_n u_{n,k_0}[x_{m_n} - e_{m_n} x_{m_n} e_{m_n}]
( {\bf 1}-p_{k_0})q_n \right\Vert_{E(\cal{M},\tau)}.
\end{align*}
We recall that $e_{k_0} \leq  p_{k_0}$ and since $\{e_n\}_{n=1}^\infty$
 is decreasing,  for $m_n
\geq k_0$, $e_{m_n} \leq p_{k_0}$ and therefore $e_{m_n} ({\bf 1}-p_{k_0}) = 0$
and since $\left\Vert q_n u_{n,k_0} \right\Vert_{\infty} \leq 1$,
we obtain that for $n$ large enough,
\begin{align*}
\left\Vert q_n \vert v_{m_n} ({\bf 1}-p_{k_0}) \vert q_n
\right\Vert_{E(\cal{M},\tau)} &= \left\Vert q_n
u_{n,k_0} (x_{m_n}) ({\bf 1}-p_{k_0}) q_n \right\Vert_{E(\cal{M},\tau)} \cr
&\leq \left\Vert x_{m_n} ({\bf 1}-p_{k_0}) q_n \right\Vert_{E(\cal{M},\tau)}.
\end{align*}
  Using Proposition~\ref{inequality},
with $x = x_{m_n}$ and $y = ({\bf 1}-p_{k_0})q_n$, we get
\begin{align*}
\left\Vert q_n \vert v_{m_n}({\bf 1}-p_{k_0}) \vert p_n
\right\Vert_{E(\cal{M},\tau)} &\leq \left\Vert
x_{m_n} \right\Vert^{\frac1{2}}_{E(\cal{M},\tau)}\
\cdot \ \left\Vert q_n (1-p_{k_0}) \vert x_{m_n} \vert
({\bf 1}-p_{k_0}) q_n \right\Vert^{\frac1{2}}_{E(\cal{M},\tau)} \cr
&\leq \left\Vert q_n ({\bf 1}-p_{k_0}) \vert x_{m_n} \vert ({\bf 1}-p_{k_0}) q_n
\right\Vert^{\frac1{2}}_{E(\cal{M},\tau)}.
\end{align*}
This implies that
$$
\limsup_{n \to \infty}
\left\Vert q_n ({\bf 1}-p_{k_0}) \vert x_{m_n} \vert
({\bf 1}-p_{k_0}) q_n \right\Vert_{E(\cal{M},\tau)} > 0.
$$
Let $s_n$ be the left support projection of
$({\bf 1}-p_{k_0})q_n$ (this is equal to the  right
support projection of $q_n ({\bf 1}-p_{k_0})$).  We have
\begin{align*}
\left\Vert q_n ({\bf 1}-p_{k_0}) \vert x_{m_n} \vert
({\bf 1}-p_{k_0})q_n \right\Vert_{E(\cal{M},\tau)} &=
\left\Vert q_n ({\bf 1}-p_{k_0}) s_n \vert x_{m_n} \vert s_n
({\bf 1}-p_{k_0}) q_n
\right\Vert_{E(\cal{M},\tau)} \cr
&\leq \left\Vert s_n \vert x_{m_n} \vert s_n \right\Vert_{E(\cal{M},\tau)}.
\end{align*}
By the definition of support projection,
$s_n \leq ({\bf 1}-p_{k_0})$ for every $n
\geq 1$, so $\{s_n\}_{n=1}^\infty$ is a sequence of finite projections.
   As in proof of Lemma~\ref{finitecompact1}, we note that
 $s_n = q_n \vee p_{k_0} - p_{k_0}$
  and as before,
$s_n \sim q_n -
q_n \wedge p_{k_0}$.
 Now since $q_n \downarrow_n 0$, $q_n -  q_n \wedge
p_{k_0}\downarrow_n 0$ hence
$\tau(s_n) = \tau(q_n - q_n \wedge p_{k_0})$
converges to zero which implies that
$s_n \downarrow_n 0$. Therefore, $\{s_n\}_{n=1}^\infty \in  \script{D}_1$. 

 In
summary, we get $\{s_n\}_{n=1}^\infty \in \script{D}_1$ with  $s_n
\leq 1-p_{k_0}$ for each $n\geq 1$ and for some $\gamma >0$,
\begin{equation}\label{gamma}
\limsup_{n \to\infty}
\left\Vert s_n \vert x_{m_n} \vert
s_n \right\Vert = \gamma.
\end{equation}
Let $f_n: = s_n \vee e_{m_n}$.

\noindent
For each $m_n \geq k_0$, $s_n \leq {\bf 1} - p_{k_0} \leq
 {\bf 1} - e_{k_0}$ so $s_n
\perp e_{m_n}$ hence $f_n = s_n + e_{m_n}$.  In
particular $\{f_n\}_{n=1}^\infty \in \script{D}_1$.

Using Proposition~\ref{lower-estimate} (it applies since $\tau(f_n) < \infty)$,
\begin{align*}
\left\Vert f_n \vert x_{m_n} \vert f_n
\right\Vert^{\frac{2q}{\alpha}}_{E(\cal{M},\tau)} &\geq
\left\Vert s_n \vert x_{m_n} \vert s_n
\right\Vert^{\frac{2q}{\alpha}}_{E(\cal{M},\tau)} +
 \left\Vert e_{m_n}
\vert x_{m_n} \vert s_n
\right\Vert^{\frac{2q}{\alpha}}_{E(\cal{M},\tau)} \cr
&+ \left\Vert s_n \vert x_{m_n} \vert e_{m_n}
\right\Vert^{\frac{2q}{\alpha}}_{E(\cal{M},\tau)} +
\left\Vert e_{m_n} \vert x_{m_n} \vert e_{m_n}
\right\Vert^{\frac{2q}{\alpha}}_{E(\cal{M},\tau)} \cr
&\geq \left\Vert s_n \vert x_{m_n} \vert s_n
\right\Vert^{\frac{2q}{\alpha}}_{E(\cal{M},\tau)} + 
\left\Vert e_{m_n} \vert x_{m_n} \vert e_{m_n}
\right\Vert^{\frac{2q}{\alpha}}_{E(\cal{M},\tau)}.
\end{align*}

Taking the limit as $n$ tends to $\infty$, one gets from (\ref{delta}) and
(\ref{gamma}) that
$\delta^{\frac{2q}{\alpha}} \geq \gamma^{\frac{2q}{\alpha}} +
\delta^{\frac{2q}{\alpha}}$.
This is   a contradiction since $\gamma > 0$. 

We conclude  that for every
$n \geq 1$, the set  $\vert V({\bf 1}-p_n) \vert$ is an 
$E$-equi-integrable set.

For the case of $\vert ({\bf 1}-p_n)V \vert$, it is enough to repeat 
the argument above for
$V^*({\bf 1}-p_n)$ using the definition of $\delta^*$ 
(instead of $\delta$). Details are left to the reader. This 
ends the proof of the 
lemma.
\end{proof}

\smallskip
We will proceed to  the proof of Theorem~\ref{main}.
  Consider two cases.

\noindent  {\it Case 1:\  Assume that $V$ is 
$E$-equi-integrable.}  

\noindent
It is
enough to set $y_n = x_{m_n} - e_{m_n} x_{m_n} e_{m_n}$ and
 $z_n = e_{m_n} x_{m_n} e_{m_n}$.

\noindent  {\it Case 2:\ Assume that $V$ is not
$E$-equi-integrable.}

\noindent
 Proposition~\ref{equi-integrability2} and 
Lemma~\ref{finite-equi-integrable} imply
that there exists $\nu>0$ such that
$$\lim\limits_{n
\to \infty} \sup_{v \in V}
\left\Vert p_n v p_n \right\Vert_{E(\cal{M},\tau)} = \nu > 0.$$
Choose a subsequence $\{v_{n_k}\}_{k=1}^\infty$ such that
\begin{equation}\label{nu} 
\lim\limits_{k \to \infty}
\left\Vert p_k v_{n_k} p_k \right\Vert_{E(\cal{M},\tau)}
 = \nu > 0.
\end{equation}

\noindent
 For each $k \geq 1$, let $w_k:
= v_{n_k} - p_k v_{n_k} p_k$ and set
$$W := \left\{w_k: \ k \geq 1 \right\}.$$
\begin{lemma}
 The set $W$ is  $E$-equi-integrable.
\end{lemma}
\begin{proof}
We note first that if $k \geq n$, then
$({\bf 1}-p_n)w_k = ({\bf 1}-p_n)v_k$ and
$w_k(1-p_n) = v_k (1-p_n)$  so for fixed $n\geq 1$,
$({\bf 1}-p_n)W = \left\{({\bf 1}-p_n)w_k:\  k <n \right\} \cup
\left\{({\bf 1}-p_n)v_k:\  k \geq n \right\}$.
Similarly,
$W({\bf 1}-p_n) = \left\{w_k({\bf 1}-p_n):\  k<n \right\} 
\cup \left\{v_k({\bf 1}-p_n):\  k \geq n \right\}.$

 Lemma~\ref{finite-equi-integrable} implies that for every $n \geq 1$,
both $\vert W({\bf 1}-p_n) \vert$ and $({\bf 1}-p_n)W$ are
$E$-equi-integrable sets. Therefore, if $W$ is not $E$-equi-integrable, there
would be a subsequence $\{w_{k(j)}\}_{j=1}^\infty$ of 
$\{w_k\}_{k=1}^\infty$ and
$\epsilon > 0$
such that
\begin{equation}\label{epsilon}
\lim_{j \to \infty} \left\Vert p_j w_{k(j)}p_j \right\Vert_{E(\cal{M},\tau
)} = \epsilon.
\end{equation}

Using Proposition~\ref{inequality} on $v_{n_{k(j)}}$ and $p_j = (p_j - p_{k(j)}) +
p_{k(j)}$,
we obtain:
\begin{align*}
\left\Vert p_j v_{n_{k(j)}} p_j \right\Vert^{\frac{2p}{\alpha}}_{E(\cal{M},\tau)}
 \geq &\left\Vert
(p_j - p_{k(j)}) v_{n_{k(j)}} (p_j -p_{k(j)})
\right\Vert^{\frac{2q}{\alpha}}_{E(\cal{M},\tau)} +
\left\Vert p_{k(j)} v_{n_{k(j)}} (p_j - p_{k(j)})
\right\Vert^{\frac{2q}{\alpha}}_{E(\cal{M},\tau)} +  \cr
&\left\Vert (p_j - p_{k(j)}) v_{n_{k(j)}} p_{k(j)}
 \right\Vert^{\frac{2q}{\alpha}}_{E(\cal{M},\tau)} +
\left\Vert (p_j - p_{k(j)}) v_{n_{k(j)}} (p_j - p_{k(j)})
\right\Vert^{\frac{2q}{\alpha}}_{E(\cal{M},\tau)}.
\end{align*}

Taking into account the identities,
$(p_j - p_{k(j)}) v_{n_{k(j)}} (p_j - p_{k(j)}) = (p_j -
p_{k(j)}) w_{k(j)} (p_j - p_{k(j)})$,
$(p_j - p_{k(j)}) v_{n_{k(j)}} p_{k(j)} = p_j w_{k(j)}
p_{k(j)}$ and
$p_{k(j)} v_{n_{k(j)}} (p_j - p_{k(j)}) = p_{k(j)} w_{k(j)}
p_j$,
one can deduce that,
\begin{align*}
\left\Vert p_j v_{n_{k(j)}} p_j
\right\Vert^{\frac{2q}{\alpha}}_{E(\cal{M},\tau)}
\geq
&\left\Vert(p_j - p_{k(j)}) w_{k(j)} (p_j -p_{k(j)})
\right\Vert^{\frac{2q}{\alpha}}_{E(\cal{M},\tau)} +
\left\Vert p_j w_{k(j)} p_{k(j)}
\right\Vert^{\frac{2q}{\alpha}}_{E(\cal{M},\tau)} + \cr
&\left\Vert p_{k(j)} w_{k(j)} p_j
\right\Vert^{\frac{2q}{\alpha}}_{E(\cal{M},\tau)} +
\left\Vert p_{k(j)} v_{n_{k(j)}} p_{k(j)}
\right\Vert^{\frac{2q}{\alpha}}_{E(\cal{M},\tau)}.
\end{align*}

Let $C(q,\alpha)$ be the norm of the identity map from
 $\ell^{\frac{2q}{\alpha}}_3$ onto
 $\ell^{\alpha}_3$, where $\ell^{\frac{2q}{\alpha}}_3$
(resp. $\ell^{\alpha}_3$) denotes the $3$-dimensional
$\ell^{\frac{2q}{\alpha}}$-space (resp. $\ell^{\alpha}$-space). We have

\begin{align*}
\left\Vert p_j v_{n_{k(j)}} p_j \right\Vert^{\frac{2q}{\alpha}}_{E(\cal{M},\tau)}
 \geq &C(q,
\alpha)^{\frac{2q}{\alpha}} \bigg[ \bigg( \left\Vert (p_j - p_{k(j)}) w_{k(j)}
(p_j - p_{k(j)}) \right\Vert^{\alpha}_{E(\cal{M},\tau)}  + \cr
 &\left\Vert p_j w_{k(j)} p_{k(j)} \right\Vert^{\alpha}_{E(\cal{M},\tau)}
  + \left\Vert p_{k(j)}
w_{k(j)} p_j \right\Vert^{\alpha}_{E(\cal{M},\tau)} \bigg)^{\frac1{\alpha}}
\bigg]^{\frac{2q}{\alpha}} + \cr
 &\left\Vert p_{k(j)} v_{n_{k(j)}} p_{k(j)}
\right\Vert^{\frac{2q}{\alpha}}_{E(\cal{M},\tau)}.
\end{align*}
We remark that $p_j w_{k(j)} p_j = (p_j - p_{k(j)}) w_{k(j)} (p_j
- p_{k(j)}) + p_j w_{k(j)} p_{k(j)} + p_{k(j)} w_{k(j)} p_j$
and since $E(\cal{M},\tau)$ is $\alpha$-convex
 (with constant $1$),
 the above inequality implies
$$\left\Vert p_j v_{n_{k(j)}} p_j \right\Vert^{\frac{2q}{\alpha}}_{E(\cal{M},\tau)}
 \geq
C(q,\alpha)^{\frac{2q}{\alpha}} \left\Vert p_j w_{k(j)} p_j
\right\Vert^{\frac{2q}{\alpha}}_{E(\cal{M},\tau)}
+ \left\Vert p_{k(j)} v_{n_{k(j)}} p_{
k(j)}
 \right\Vert^{\frac{2q}{\alpha}}_{E(\cal{M},\tau)},
$$
and taking the limit as $j \to \infty$, we get from (\ref{nu}) and (\ref{epsilon})
that
$$
\nu^{\frac{2q}{\alpha}} \geq C(q, \alpha)^{\frac{2q}{\alpha}}
\epsilon^{\frac{2q}{\alpha}} + \nu^{\frac{2q}{\alpha}}.
$$
This is a contradiction since
$\epsilon >0$, so $W$ is a $E$-equi-integrable set.
The lemma is proved.

\smallskip

To complete the proof of Theorem~\ref{main}, we note  that
$W = \{x_{n_{k}} - p_k x_{n_{k}} p_k: \ k \geq 1\}$ so if we set $y_k =
x_{n_{k}} - p_k x_{n_{k}} p_k$ and $z_k = p_k x_{n_{k}} p_k$. The proof of
 is complete.
\end{proof}

\begin{remarks}
(1) If $\M =\cal{B}(\ell^2)$ with the usual trace, then every
projection of finite trace is a finite rank projection so 
in the proof above, $\delta=\delta^*=0$. In the particular case
of unitary matrix space $C_E$ where  $E$ is a symmetric 
sequence space, one proceed directly to Case~2 by setting
$W:=\left\{ x_n - p_n x_n p_n;\ n \geq 1 \right\}$ where
$\{p_n\}_{n=1}^\infty$ is an arbitrary sequence of projections satisfying:
$p_n \downarrow_n 0$ and for every $n \geq 1$, $1-p_n$ is a finite
projection.

(2) If $\M$ is a finite von Neumann algebra with 
a normalized finite trace $\tau$ and $E$ is a symmetric space
on $[0,1]$ satisfying the assumptions of Theorem~\ref{main}, it is enough to take
$p_n =e_n$ (i.e  $g_n ={\bf 1}-e_n$ on Lemma~3.2) and conclude 
immediately as in Lemma~\ref{finite-equi-integrable} that $V$ is $E$-equi-integrable.

 (3) In the proof above, it is clear that
the projections $\{p_k\}_{k=1}^\infty$ are such that
either $\T(p_1)< \infty$ or $\tau({\bf 1}-p_k)< \infty$ for all $k
\geq 1$. In fact, the argument above shows that if $\{e_n\}_{n=1}^\infty$ is a sequence 
in $\cal{D}_1$ that attained the quantities $\delta$ and $\delta^*$, then any
sequence of projections satisfying $p_n \downarrow_n 0$, $e_n \leq p_n$ for
each $n\geq 1$ and $\T({\bf 1}-p_n)<\infty$ for each $n\geq 1$, would
satisfy the conclusion of Theorem~\ref{main}.
\end{remarks}

The following extension shows that if one considers finitely 
many bounded sequences in $E(\M,\T)$, one
can choose a single sequence of projections that works for
each sequence. 

\begin{corollary}
If $\M$ and $E$ are as in Theorem~\ref{main} and 
$\left\{x^{(1)}\right\}_{n=1}^\infty,\left\{x^{(2)}\right\}_{n=1}^\infty, 
\dots, 
\left\{x^{(j_0)}\right\}_{n=1}^\infty$ be finitely many bounded sequences in
$E(\M, \tau)$. Then there exist a strictly increasing sequence 
$\{n_k\}_{k=1}^\infty$ of $\N$ 
 and a sequence of decreasing projections
 $p_k \downarrow_k 0$ in $\cal M$ such that  for
each $1 \leq j \leq j_0$, the set
 $\left\{x^{(j)}_{n_k}- p_k x^{(j)}_{n_k}p_k :\ k \geq 1 \right\}$ 
is 
$E$-equi-integrable.
\end{corollary}
\begin{proof}
  For   $1 \leq j \leq j_0$, we set, 
as in the proof of Theorem~\ref{main}, 
 $$\delta_j: = \sup\left\{\lim\limits_{n \to \infty} 
\sup\limits_{k \in \N}
\left\Vert e_n \vert x^{(j)}_k \vert e_n \right\Vert_{E(\cal{M},\tau)}:\ 
\{e_n\}_{n=1}^\infty \in \script{D}_1 \right\}.$$

 One can choose a strictly increasing sequence $\{n_k\}_{k=1}^\infty$ 
in $\N$ such that for each 
$1\leq j \leq j_0$, there exists a sequence 
  $\{e^{(j)}_k\}_{k=1}^\infty \in \script{D}_1$ 
with
$$\delta_j =\lim_{k \to \infty} \left\Vert e^{(j)}_k
\vert x^{(j)}_{n_k} \vert e^{(j)}_k
\right\Vert_{E(\cal{M},\tau)}$$
 and
$$\delta^*_j=\lim\limits_{k \to \infty} \left\Vert e^{(j)}_k 
\vert {x^{(j)}_{n_k}}* \vert e^{(j)}_k
\right\Vert_{E(\cal{M},\tau)}
= \sup \left\{\lim\limits_{n \to \infty} \sup\limits_{k} 
\left\Vert q_n \vert
{x^{(j)}_{n_k}}^* \vert q_n \right\Vert_{E(\cal{M},\tau)}; 
\{q_{n}\}_{n=1}^\infty \in \script{D}_1 \right\}.$$
For every $k \geq 1$, set 
$e_k := \vee_{1\leq j \leq j_0} e^{(j)}_k$. Since 
$\T(e_k) \leq \sum_{j=1}^{j_0} \T(e^{(j)}_k)$, it is clear that 
the sequence $\{e_k\}_{k=1}^\infty$ belongs to $\script{D}_1$ and each of the 
$\delta_j$'s and $\delta^{*}_j$'s are attained at $\{e_k\}_{k=1}^\infty$.
One can complete the proof by procceding as in the proof
of Theorem~\ref{main}, simmultaneously on the finite set of sequences and
the fixed $\{e_k\}_{k=1}^\infty$. 
\end{proof}

\bigskip

Our next result shows that the decreasing projections in the decomposition
can be replaced by mutually disjoint projections.

\begin{theorem}\label{main2}
Let $E$ be an order continuous quasi-Banach function 
space as in Theorem~\ref{main}.
  Let $\{x_n\}_{n=1}^\infty$ be a bounded sequence in $E(\cal{M},\tau)$ then there exists a
subsequence $\{x_{n_{k}}\}_{k=1}^\infty$ of $\{x_n\}_{n=1}^\infty$, 
bounded sequences $\{\varphi_k\}_{k=1}^\infty$ and
$\{\zeta_k\}_{k=1}^\infty$ in $E(\cal{M,\tau})$ and mutually disjoint sequence of projections
$\{e_k\}_{k=1}^\infty$ such that:

\begin{itemize}
\item [(i)] $x_{n_{k}} = \varphi_k + \zeta_k$  for all $k \geq 1$;
\item [(ii)] $\{\varphi_k:\  k\geq 1\}$ is $E$-equi-integrable and
$e_k \varphi_k e_k = 0 $ for all $\ k \geq 1$;
\item [(iii)] $\{\zeta_k\}_{k=1}^\infty$ is such that $e_k \zeta_k e_k = \zeta_k$
 for all $k \geq1$.
\end{itemize}
\end{theorem}

\begin{proof}
Let $\{x_n\}_{n=1}^\infty$  be a bounded sequence in $E(\cal{M},\tau)$
and suppose (by taking a subsequence if
necessary), $x_n = y_n + z_n$ with $p_n y_n p_n = 0$, the set
$\{y_n:\  n \geq 1\}$ is  $E$-equi-integrable and $p_n z_n p_n = z_n$
for all $n \geq 1$,  be the decomposition of $\{x_n\}_{n=1}^\infty$ as  
in Theorem~\ref{main}.

Let $n_1 = 1$. Since $p_n \downarrow_n 0$ and
$$
p_1 z_1 p_1 - (p_1 - p_n) z_1 (p_1 - p_n) = p_n z_1 p_1 + p_1
 z_1 p_n - p_n z_1 p_n,
$$
Proposition~ref{order-continuity}(part(ii)) shows that
$$
\lim_{n \to \infty} \Vert p_1 z_1 p_1 - (p_1 - p_n) z_1 (p_1 -
p_n) \Vert_{E(\cal{M},\tau)} = 0.
$$
Choose $n_2 > n_1 = 1$ such that
$$
\Vert p_1 z_1 p_1 - (p_1 - p_{n_{2}}) z_1 (p_1 - p_{n_{2}})
\Vert_{E(\cal{M},\tau)} < \frac1{2}.
$$
Inductively, one can construct $n_1 < n_2 < \dots < n_k < \dots$
such that
$$\Vert p_{n_{k}} z_{n_{k}} p_{n_{k}} - (p_{n_{k}} - p_{n_{k + 1}})
z_{n_{k}} (p_{n_{k}} - p_{n_{k + 1}}) \Vert_{E(\cal{M},\tau)}
 < {\frac1{2^k}}.
$$
Since $z_n=p_n z_n p_n$  for every $n\geq 1$, one gets
$$
\Vert z_{n_{k}} - (p_{n_{k}} - p_{n_{k + 1}}) z_{n_{k}} (p_{n_{k}} - p_{n_{k +
1}}) \Vert_{E(\cal{M},\tau)} < {\frac1{2^k}}.
$$
For every $ k \geq 1$, set
\begin{align*}
e_k : & = p_{n_{k}} - p_{n_{k + 1}} \cr
\zeta_k : & = (p_{n_{k}} - p_{n_{k + 1}}) z_{n_{k}} (p_{n_{k}} - p_{n_{k +
1}}) \cr
\varphi_k : & = y_{n_{k}} + [z_{n_{k}} - e_k z_{n_{k}} e_k] .
\end{align*}
Since $\{y_{n_{k}}:\ k \geq 1\}$ is a $E$-equi-integrable set and
$\lim\limits_{k \to \infty} \Vert z_{n_{k}} - e_k z_{n_{k}} e_k
\Vert_{E(\cal{M},\tau)} = 0$, it is clear that
 $\{\varphi_k:\  k \geq 1\}$ is $E$-equi-integrable.
Also $\{e_k\}_{k=1}^\infty$ is mutually disjoint.  The proof is complete.
\end{proof}

\begin{corollary}
Let $E$ be an order-continuous symmetric Banach function space 
on $\R^{+}$ with the Fatou propery.
  Let $\{x_n\}_{n=1}$ be a bounded sequence in $E(\cal{M},\tau)$ then there exists a
subsequence $\{x_{n_{k}}\}_{k=1}^\infty$ of $\{x_n\}_{n=1}^\infty$, 
bounded sequences $\{\varphi_k\}_{k=1}^\infty$ and
$\{\zeta_k\}_{k=1}^\infty$ in $E(\cal{M,\tau})$ and mutually disjoint sequence of projections
$\{e_k\}_{k=1}^\infty$ such that:
\begin{itemize}
\item [(i)] $x_{n_{k}} = \varphi_k + \zeta_k$  for all $k \geq 1$;
\item [(ii)] $\{\varphi_k:\  k\geq 1\}$ is $E$-equi-integrable and
$e_k \varphi_k e_k = 0 $ for all $\ k \geq 1$;
\item [(iii)] $\{\zeta_k\}_{k=1}^\infty$ is such that $e_k \zeta_k e_k = \zeta_k$
 for all $k \geq1$.
\end{itemize}
\end{corollary}

\begin{proof}
 Assume that   $E$ has the Fatou
property (equivalenty $E$ does not contain $c_0$). Since $E$ is 
symmetric, $E \not \supset c_0$ is equivalent to $E$ not containing
$\ell^n_\infty$ uniformly, and therefore $E$ satisfies the $q$-lower
estimate for some $q$ and one can renorm $E$ so that it satisfies the
lower $q$-estimate of constant $1$.  All of these facts can be found in
\cite{LT}.
\end{proof}

\bigskip

 The proof of Theorem~\ref{main} can be adjusted to obtain decompositions where
the projections are taken only on one side, that is,  the following
result follows:
\begin{corollary}
Let $E$ be an order continuous quasi-Banach function space in $\R^+$ that
is $\alpha$-convex with constant $1$ for some $0 < \alpha \leq 1$ and
suppose that $E$ satisfies a lower $q$-estimate with constant $1$ for some
$q \geq \alpha$.
Let $\{x_n\}_{n=1}^\infty$ be a bounded sequence in $E(\cal{M},\tau)$
 then there exist a
subsequence $\{x_{n_k}\}_{k=1}^\infty$ of $\{x_n\}_{n=1}^\infty$, bounded 
sequences $\{y_k\}_{k=1}^\infty$ and
$\{z_k\}_{k=1}^\infty$ in $E(\cal{M},\tau)$ and decreasing projections
 $e_k \downarrow_n 0$ in $\cal M$ such that:
\begin{itemize}
\item [(i)] $x_{n_k} = y_k + z_k$ for all $k \geq 1$;
\item [(ii)]  $e_ky_k = 0$
for all $k \geq 1$ and
$\lim\limits_{n\to \infty}\sup_{k\geq 1}||f_n y_k||_{E(\M,\T)}=0$ for
every $f_n \downarrow_n 0$.
\item [(iii)] $\{z_k\}_{k=1}^\infty$ is such that $e_kz_k = z_k$ for all $k \geq
1$.
\end{itemize}
\end{corollary}

\begin{definition}
A subspace $X$  of $L^p(\M,\T)$ is  called strongly 
embedded into $L^p(\M,\T)$ if the $L^p$ and the measure topologies
on $X$ coincide.
\end{definition}

The following result is a direct application of 
Proposition~\ref{measure} and Theorem~\ref{main2}.
\begin{theorem}
Let $1\leq p <\infty$. Every subspace
of $L^p(\M,\T)$ either contains almost isometric
copies of $\ell^p$ or is strongly embedded in $L^p(\M,\T)$.
\end{theorem}

The next corollary should be compared with
\cite[Theorem~2.4]{SU}.

\begin{corollary} Assume that $\M$ is finite and
$p>2$. Every subspace of $L^p(\M,\T)$  either
contains almost isometric copies of $\ell^p$ or 
is isomorphic to a Hilbert space.
\end{corollary}

\medskip

For the commutative case, the space $\ell^p$ can not be
strongly embedded in $L^p[0,1]$ for $0<p<2$. This is
due to Kalton \cite{KA3} for $0<p<1$ and Rosenthal
\cite{R6} for the case $1\leq p <2$ (see also \cite{CADI} for
another approach). A  non-commutative analogue should be
of interest.

\medskip

\noindent
{\bf Problem:}
{\it  Let $\M$ be a semifinite von Neumann algebra and
$0<p<2$.
Does $\ell^p$ strongly embed into $L^p(\M,\T)$?}


\providecommand{\bysame}{\leavevmode\hbox to3em{\hrulefill}\thinspace}

\end{document}